\pdfoutput=1 

\documentclass[oneside,english,12pt,titlepage]{amsbook}
\usepackage[T1]{fontenc}
\usepackage[utf8]{luainputenc}
\usepackage{verbatim}
\usepackage{amsthm}
\usepackage{amstext}
\usepackage{amssymb}
\usepackage{graphicx}

\usepackage{amscd}

\makeatletter

\let\SF@@footnote\footnote
\def\footnote{\ifx\protect\@typeset@protect
    \expandafter\SF@@footnote
  \else
    \expandafter\SF@gobble@opt
  \fi
}
\expandafter\def\csname SF@gobble@opt \endcsname{\@ifnextchar[
  \SF@gobble@twobracket
  \@gobble
}
\edef\SF@gobble@opt{\noexpand\protect
  \expandafter\noexpand\csname SF@gobble@opt \endcsname}
\def\SF@gobble@twobracket[#1]#2{}
\providecommand{\tabularnewline}{\\}

\theoremstyle{plain}
\newtheorem{thm}{Theorem}[chapter]

\newtheorem{lem}[thm]{Lemma}
\newtheorem{defn}[thm]{Definition}

\newtheorem*{thm*}{Theorem}
\newtheorem*{cor*}{Corollary}
\newtheorem*{lem*}{Lemma}
\newtheorem*{defn*}{Definition}
\newtheorem*{prop*}{Proposition}
\newtheorem*{conj*}{Conjecture}

\makeatother

\usepackage{babel}

\begin{document}

\pagestyle{plain}

\frontmatter

\begin{titlepage}

{\Huge The 18-cycle in Bianchi $VI_{-1/9}^{^{*}}$:}

\vspace{0.5cm}

{\Huge Combined Linear Local Passage and}

\vspace{0.5cm}

\begin{center}

{\Huge Numerical Simulation}

\vspace{0.5cm}


\vspace{1cm}

{\Large Johannes Buchner} 

\vspace{0.5cm}

\today

\vspace{2cm}

{\bf  Abstract}

\end{center}

\vspace{0.5cm}

In this paper, we find an example for a periodic heteroclinic
chain in Bianchi $VI_{-1/9}^{^{*}}$ that allows Takens Linearization
at all base points. It turns out to be a "18-cycle'', i.e.
involving a heteroclinic chain of 18 different base points at the
Kasner circle. We then show that the combined cinear local passage
at the 18-cycle is a contraction. This qualifies the 18-cycle as a candidate for proving the first rigorous convergence theorem in Bianchi $VI_{-1/9}^{^{*}}$.

For the numerical simulation of solutions that follow such heteroclinic cycles, we use a variable-step, variable-order (VSVO) Adams-Bashforth-Moulton PECE solver in Matlab. We conclude with a discussion on how to proceed further in studying Bianchi cosmologies, and also discuss directions for future research in  inhomogeneous (PDE-) cosmological models. This puts our results in a broader perspective. The appendix contains symbolic and computations done by Mathematica for examples discussed throughout the text.

\end{titlepage}

\newpage{}

\tableofcontents{}

\newpage{}

\mainmatter

\newpage

\chapter{Bianchi Spacetimes - Existing Results, Challenges and Techniques}
\label{bianchi-defs}

\section{The Equations of Wainwright and Hsu}
In the paper \cite{waihsu89} by Wainwright and Hsu, a formulation of the Einstein Equations for Bianchi models is presented that has several advantages, one of them beeing that it contains all models of Bianchi class A. Here are the equations, which are used throughout this dissertation:

\medskip

\begin{equation}\label{bix-first}
\begin{array}{rcl}
   N_1'      &=& ( q - 4 \Sigma_+                     ) N_1,
\\ N_2'      &=& ( q + 2 \Sigma_+ + 2\sqrt{3}\Sigma_- ) N_2,
\\ N_3'      &=& ( q + 2 \Sigma_+ - 2\sqrt{3}\Sigma_- ) N_3,
\\ \Sigma_+' &=& (q-2) \Sigma_+ - 3 S_+,
\\ \Sigma_-' &=& (q-2) \Sigma_- - 3 S_-.
\end{array}
\end{equation}
%

with constraint

\begin{equation}
\label{bix-constraint}
\Omega    = 1 - \Sigma_+^2 - \Sigma_-^2 -K  
\end{equation}


and abbreviations
%

\begin{equation}\label{bix-last}
\begin{array}{rcl}
   q         &=& 2 \left( \Sigma_+^2 + \Sigma_-^2 \right)
                 + \frac{1}{2}(3\gamma-2)\Omega,
\\ K         &=& \frac{3}{4} \left( \rule[2ex]{0pt}{0pt} N_1^2 + N_2^2 + N_3^2
                                    - 2\left( N_1 N_2 + N_2 N_3 + N_3 N_1 \right) 
                             \right),
\\ S_+       &=& \frac{1}{2}\left( 
                   \left( N_2 - N_3 \right)^2
                   - N_1 \left( 2 N_1 - N_2 - N_3 \right) 
                 \right),
\\ S_-       &=& \frac{1}{2}\sqrt{3} 
                 \left( N_3 - N_2 \right) \left( N_1 - N_2 - N_3 \right).
\end{array}
\end{equation}

\smallskip

The fixed parameter $\gamma$  is related to the choice of matter model (e.g. $\gamma=1$ represents dust, whereas $\gamma=4/3$ represents radiation).

The properties of equations above have been studied intensively, for a recent survey see \cite{heiugg09b}. The main goal are rigorous results on the correspondence of iterations of the so-called "Kasner map" $f$ to the dynamics of nearby trajectories to the Bianchi system (\ref{bix-first}) with reversed time, i.e.~in the $\alpha$-limit $t\to -\infty$. We will introduce the necessary background in the rest of this section.

\newpage

\subsection{Vacuum Models of Bianchi Class A}
\label{vac-a}

From now on, we will restrict ourselves to the vacuum case $\Omega=0$ . This yields a 4-dimensional model, as we have five variables and one constraint:

\begin{defn} (Phase Space of the Vacuum Wainwright-Hsu ODEs)
\[
\mathcal{W}=\{(N_{1},N_{2},N_{3},\Sigma_{+},\Sigma_{-}) \; | \; 0=1 - \Sigma_+^2 - \Sigma_-^2 -K\}
\]
\end{defn} 

As we are interested in the dynamics of the Bianchi system (\ref{bix-first}) with reversed time, i.e.~in the $\alpha$-limit $t\to -\infty$, we will we denote by $X^{W}$ the vector field corresponding to this time direction, for use in later chapters\footnote{in the paper by Béguin \cite{beg10}, the equations are presented directly with time direction chosen towards the big bang, but we stick to the form of the equations used in the classic reference \cite{waihsu89}, and also in \cite{rin01, lieetal10}.}. This means $X^{W}$ stands for the vector field corresponding to the right side of (\ref{bix-first}), multiplied by $-1$.

When we look at those equations, we observe that if $N_1=N_2=N_3=0$, the vector field is zero, as $K=0$ and $q=2$ in this case. We denote by $\mathcal{K} = \{ N_1=N_2=N_3=0,\; \Omega=0 \}$ the resulting circle of equilibria: we obtain a circle because the constraint (\ref{bix-constraint}) reduces to $\Sigma_+^2 + \Sigma_-^2=1$. It is called the "Kasner circle", because the points $p \in \mathcal{K}$ represent the Kasner solution of the Einstein Equations discussed in section \ref{kasner-solution}.

In the classification of spatially homogenous models (based on the classification of 3-dimensional Lie-Algebras by Bianchi \cite{bianchi-orginal}), these are of Bianchi class I. One advantage of the Wainwright-Hsu equations is that they contain all models of Bianchi class A, with the signs of the $N_{i}$ determining the type of Bianchi model, see the table above.

\begin{table}
\begin{tabular}{|c|ccc|}
\hline
Bianchi Class & $N_1$ & $N_2$ & $N_3$ \\
\hline
I       & $0$&$0$&$0$ \\
II      & $+$&$0$&$0$ \\
VI$_0$  & $0$&$+$&$-$ \\
VII$_0$ & $0$&$+$&$+$ \\
VIII    & $-$&$+$&$+$ \\
IX      & $+$&$+$&$+$ \\
\hline
\end{tabular}
\end{table}


If we allow one of the $N_i$ to be non-zero, the resulting half ellipsoids are called "Kasner caps": $C_k = \{ N_k > 0, \; N_l=N_m=0, \; \Omega=0 \}$ with $\{k,l,m\}=\{1,2,3\}$. They consist of heteroclinic orbits to equilibria on the Kasner circle and are of Bianchi class II. The projections of the trajectories of Bianchi class-II vacuum solutions onto the $\Sigma_\pm$-plane yield straight lines connecting two points of the Kasner circle. 

These can be constructed geometrically in the following way: for a point $p \in \mathcal{K}$ of the Kasner circle, identify the nearest corner of the circumscribed triangle, and draw the resulting line as illustrated in the picture below:

\medskip
\begin{minipage}{0.7\textwidth}
\centering
\setlength{\unitlength}{0.277\textwidth}
\begin{picture}(3.6,3.6)(-1.3,-1.8)
\put(-1.3 , 1.8 ){\makebox(0,0)[tl]{%
  \includegraphics[width=3.6\unitlength]{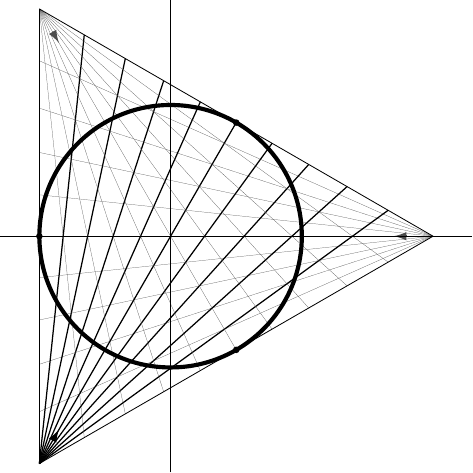}}}
\put( 2.3 ,-0.03){\makebox(0,0)[rt]{$\Sigma_+$}}
\put( 0.03, 1.8){\makebox(0,0)[lt]{$\Sigma_-$}}

\put(-1.03,-0.03){\makebox(0,0)[rt]{$T_1$}}
\put( 0.53, 0.88){\makebox(0,0)[lb]{$T_2$}}
\put( 0.53,-0.88){\makebox(0,0)[lt]{$T_3$}}

\end{picture}
\end{minipage}
\medskip

Observe that this works for any $p\in\mathcal{K}$ except for the three points where the Kasner circle touches the circumscribed triangle, which we denote by $T_1,T_2,T_3$ and refer to them as Taub points. 

\subsection{The Kasner Map}
\label{kasner-map}

This leads to the definition of the so-called Kasner map $f:\mathcal{K}\to\mathcal{K}$.
For each point $p_+\in (\mathcal{K}\setminus\{T_1,T_2,T_3\})$ there exists 
a Bianchi class-II vacuum heteroclinic orbit $H(t)$ converging to $p_+$ as $t\to\infty$. 
This orbit is unique up to reflection $(N_1,N_2,N_3)\mapsto(-N_1,-N_2,-N_3)$. 
Its unique $\alpha$-limit $p_-$ defines the image of $p_+$ under the Kasner map
\begin{equation}\label{eqKasnerMap}
  f(p_+) := p_-
\end{equation}
Including the three fixed points, $f(T_k):=T_k$, this construction yields a continuous map, $f:\mathcal{K}\to\mathcal{K}$. In fact $f$ is a non-uniformly expanding map and its image $f(\mathcal{K})$ is a double cover of $\mathcal{K}$, which can be seen directly from the geometric description given above.

For use in later chapters, let us denote by $H_{p,f(q)}$ the heteroclinic Bianchi-II-orbit from $p$ to its image point under the Kasner map $f(q)$ and by $H_B=\bigcup_{p\in B}H_{p,f(q)}$ the set of all heterclinic Bianchi-II-orbits connecting two basepoints in the set $B \subset \mathcal{K}$.

\newpage 

We will now introduce the so-called Kasner-parameter $u$, as it is a convenient way to parametrize the Kasner circle. Let us divide the Kasner circle into six sectors and label them as follows:

\smallskip
\begin{minipage}{0.7\textwidth}
\centering
\setlength{\unitlength}{0.277\textwidth}
\begin{picture}(3.6,3.6)(-1.3,-1.8)
\put(-1.3 , 1.8 ){\makebox(0,0)[tl]{%
  \includegraphics[width=3.6\unitlength]{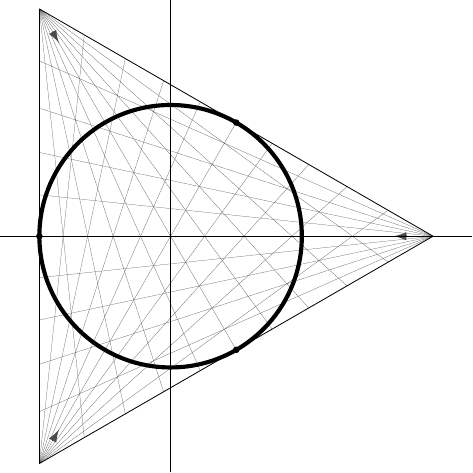}}}
\put( 2.3 ,-0.03){\makebox(0,0)[rt]{$\Sigma_+$}}
\put( 0.03, 1.8){\makebox(0,0)[lt]{$\Sigma_-$}}

\put(-1.03,-0.03){\makebox(0,0)[rt]{$T_1$}}
\put( 0.53, 0.88){\makebox(0,0)[lb]{$T_2$}}
\put( 0.53,-0.88){\makebox(0,0)[lt]{$T_3$}}

\put(1.03,-0.03){\makebox(0,0)[lt]{$Q_1$}}
\put( -0.53, 0.88){\makebox(0,0)[rb]{$Q_2$}}
\put( -0.53,-0.88){\makebox(0,0)[rt]{$Q_3$}}


\put( 0,1.03){\makebox(0,0)[lb]{{\bf \Large 2 }}}
\put( 0,-1.03){\makebox(0,0)[lt]{{\bf \Large 5 }}}
\put( 0.88, 0.53){\makebox(0,0)[bl]{{\bf \Large 1 }}}
\put( -0.88,0.53){\makebox(0,0)[rb]{{\bf \Large 3 }}}
\put( -0.88, -0.53){\makebox(0,0)[tr]{{\bf \Large 4}}}
\put( 0.88,-0.53){\makebox(0,0)[lt]{{\bf \Large 6 }}}

\end{picture}
\end{minipage}
\smallskip

The Kasner parameter ranges in $[1,\infty]$, where it holds that $u=\infty$ at the Taub point $T_i$ and $u=1$ at the points $Q_i$ shown in the picture. So for each value of $u$, we get an equivalence class of six points cooresponding to the $p \in \mathcal{K}$ in each sector of the Kasner circle (except for the values $u=\infty$ and $u=1$, where the equivalence class consists only of three points). Expressed in $u$, the Kasner map has a very simple form:

\begin{eqnarray*}
f(u) & = & \begin{cases}
u-1 & u\in[2,\infty]\\
\frac{1}{u-1} & u\in[1,2]
\end{cases}
\end{eqnarray*}


In the picture below, the dynamics of the Kasner-map is shown: If you start close to a Taub point (meaning that u is large compared to 1), there are first "bounces" around this Taub point as the value of u is decreased by 1 in each step. Then, after the value of u has fallen below 2, there is a "excursion" to a different part of the Kasner circle: 

\bigskip
\begin{minipage}{0.7\textwidth}
\centering
\setlength{\unitlength}{0.277\textwidth}
\begin{picture}(3.6,3.6)(-1.3,-1.8)
\put(-1.3 , 1.8 ){\makebox(0,0)[tl]{%
  \includegraphics[width=3.6\unitlength]{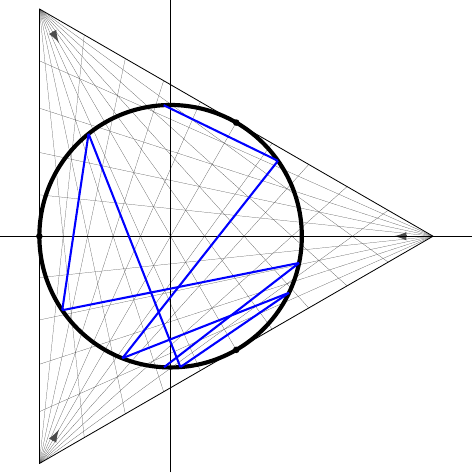}}}
\put( 2.3 ,-0.03){\makebox(0,0)[rt]{$\Sigma_+$}}
\put( 0.03, 1.8){\makebox(0,0)[lt]{$\Sigma_-$}}

\put(-1.03,-0.03){\makebox(0,0)[rt]{$T_1$}}
\put( 0.53, 0.88){\makebox(0,0)[lb]{$T_2$}}
\put( 0.53,-0.88){\makebox(0,0)[lt]{$T_3$}}

\end{picture}
\end{minipage}
\bigskip

There is also an interesting connection to the Kasner solution described in section \ref{kasner-solution}. Each sector of the Kasner circle corresponds to a permutation of the Kasner expontents $p_{i}$, which can be expressed in the Kasner parameter $u$:

\bigskip
\begin{minipage}{0.7\textwidth}
\centering
\setlength{\unitlength}{0.277\textwidth}
\begin{picture}(3.6,3.6)(-1.3,-1.8)
\put(-1.3 , 1.8 ){\makebox(0,0)[tl]{%
  \includegraphics[width=3.6\unitlength]{pics/BIX-only-circle}}}
\put( 2.3 ,-0.03){\makebox(0,0)[rt]{$\Sigma_+$}}
\put( 0.03, 1.8){\makebox(0,0)[lt]{$\Sigma_-$}}

\put(-1.03,-0.03){\makebox(0,0)[rt]{$T_1$}}
\put( 0.53, 0.88){\makebox(0,0)[lb]{$T_2$}}
\put( 0.53,-0.88){\makebox(0,0)[lt]{$T_3$}}

\put(1.03,-0.03){\makebox(0,0)[lt]{$Q_1$}}
\put( -0.53, 0.88){\makebox(0,0)[rb]{$Q_2$}}
\put( -0.53,-0.88){\makebox(0,0)[rt]{$Q_3$}}

\put( 0,1.03){\makebox(0,0)[b]{{ \Large (213)}}}
\put( 0,-1.03){\makebox(0,0)[t]{{ \Large (312)}}}

\put( 0.88, 0.53){\makebox(0,0)[bl]{{ \Large (123)}}}
\put( -0.88,0.53){\makebox(0,0)[rb]{{ \Large (231)}}}
\put( -0.88, -0.53){\makebox(0,0)[tr]{{ \Large (321)}}}
\put( 0.88,-0.53){\makebox(0,0)[lt]{{ \Large (132)}}}



\end{picture}
\end{minipage}
\bigskip

 where sector (321) means e.g. that $p_{3}<p_{2}<p_{1}$ which fixes the formula for each of them. As an example, consider the sector 5 or (312) (for details see \cite{heiugg09b}, p.8): 

\bigskip

\[
p_{3}=\frac{-u}{1+u+u^{2}}
\]

\[
p_{1}=\frac{(u+1)}{1+u+u^{2}}
\]

\[
p_{2}=\frac{u(u+1)}{1+u+u^{2}}
\]

\bigskip
\bigskip

We will need the formulas above, as they will allow us to express the eigenvalues of the linearized vector field at points of the Kasner circle in $u$, a key step for obtaining our results in the later chapters.

\subsection{Eigenvalues in Terms of the Kasner Parameter $u$}
\label{ev-formula-section}

When we linearize the vector field corresponding to equations (\ref{bix-first}) at points of the Kasner circle, we arrive at the following Matrix:

\bigskip
\medskip

\begin{equation*}
\left(\begin{array}{ccccc}
2 - 4 \Sigma_+ \hspace{-1em}& 0 & 0 & 0 & 0 \\
0 & 2 + 2 \Sigma_+ + 2\sqrt{3} \Sigma_- \hspace{-2em}& 0 & 0 & 0 \\
0 & 0 & 2 + 2 \Sigma_+ - 2\sqrt{3} \Sigma_- \hspace{-2em}& 0 & 0 \\
0 & 0 & 0 & 3(2-\gamma)\Sigma_+^2 & 3(2-\gamma)\Sigma_+\Sigma_- \\
0 & 0 & 0 & 3(2-\gamma)\Sigma_+\Sigma_- & 3(2-\gamma)\Sigma_-^2
\end{array}\right).
\end{equation*}

\bigskip
\bigskip

and we can compute the following eigenvalues to eigenvectors $\partial_{N_1}$, $\partial_{N_2}$, $\partial_{N_3}$ tangential to the Bianchi class-II vacuum heteroclinics:

\[
\begin{array}{rcl}
\mu_1 &=& 2-4\Sigma_+,\\
\mu_2 &=& 2 + 2 \Sigma_+ + 2\sqrt{3}\Sigma_-, \\
\mu_3 &=& 2 + 2 \Sigma_+ - 2\sqrt{3}\Sigma_-
\end{array}
\]

\bigskip

In addition we have the trivial eigenvalue zero to the eigenvector 
$-\Sigma_-\partial_{\Sigma_+} + \Sigma_+\partial_{\Sigma_-}$ 
tangential to the Kasner circle $\mathcal{K}$.
The fifth eigenvalue $\mu_\Omega = 3(2-\gamma) > 0$ corresponds to the eigenvector 
$\Sigma_+\partial_{\Sigma_+} + \Sigma_-\partial_{\Sigma_-}$ transverse to the vacuum boundary $\{\Omega=0\}$.

\newpage

Now we use that it is possible to express the $\Sigma_{+/-}$-variables
in terms of the Kasner exponents $p_i$ (see \cite{heiugg09b}, p.7):

\[
\Sigma_{+}=-\frac{3}{2}p_{1}+\frac{1}{2}
\]

\[
\Sigma_{-}=-\frac{\sqrt{3}}{2}(p_{1}+2p_{2}-1)
\]

Thus, we arrive at the following formulas for the eigenvalues expressed in $u$:

\begin{equation}
\left(\lambda_{1},\lambda_{2},\lambda_{3}\right)=\left(\frac{-6u}{1+u+u^{2}},\frac{6(1+u)}{1+u+u^{2}},\frac{6u(1+u)}{1+u+u^{2}}\right)
\label{eq:EVinBIXinU}
\end{equation}

As it holds that $u \in [1,\infty]$, we observe that at each point of the Kasner circle, there is one negative and two positive eigenvalues. But recall that we are interested in the time-direction $t\to -\infty$. The negative eigenvalue is unstable towards the past, while the two positive eigenvalues are stable in backwards-time. This means that for our vector field $X^W$ (which has the time-direction already reversed, see beginning of section \ref{vac-a}) there is one unstable eigenvalue $\lambda_{u}$ and two stable eigenvalues $\lambda_{s}, \lambda_{ss}$, at a point of the Kasner circle. Away from the Taub points it also holds that:

\[
|\lambda_{u}|<|\lambda_{s}|<|\lambda_{ss}|
\]

\bigskip

Finally, also note that in Bianchi IX, the situation in different sectors of the Kasner-circle only differs by a permutation of those 3 formulas for the eigenvalues. This makes it easy to examine the question of resonances of the eigenvalues, which we are trying to exclude when linearizing the vector field.

However, we will also deal with Bianchi $VI_{-1/9}^{^{*}}$ later, and there the situation is more complicated as the formulas for the eigenvalues at points on the Kasner circle do depend on the sector, so in order to check for resonances, a lot of cases have to be considered. This will be done in chapter \ref{bvi-chapter}. In the next section, we will introduce the equations for Bianchi $VI_{-1/9}^{^{*}}$, which is of class B and not covered by the equations of Wainwright and Hsu.

\newpage

\section{Bianchi $VI_{-1/9}^{^{*}}$}
\label{bianchi-vi-defs}

We now present the equations for Bianchi $VI_{-1/9}^{^{*}}$, which is the most general model in Bianchi class B, and has a crucial importance for inhomogenous cosmologies (see e.g. \cite{hew03}):

\begin{align}
\Sigma_{+}^{'} & =(q-2)\Sigma_{+}+3\Sigma_{2}^{2}-2N_{-}^{2}-6A^{2}
\end{align}
 
\begin{align}
\Sigma_{-}^{'} & =(q-2)\Sigma_{-}-\sqrt{3}\Sigma_{2}^{2}+2\sqrt{3}\Sigma_{\times}^{2}-2\sqrt{3}N_{-}^{2}+2\sqrt{3}A^{2}
\end{align}

\begin{align}
\Sigma_{\times}^{'} & =(q-2-2\sqrt{3}\Sigma_{-})\Sigma_{\times}-8N_{-}A
\end{align}

\begin{align}
\Sigma_{2}^{'} & =(q-2-3\Sigma_{+}+\sqrt{3}\Sigma_{-})\Sigma_{2}
\end{align}

\begin{align}
N_{-}^{'} & =(q+2\Sigma_{+}+2\sqrt{3}\Sigma_{-})N_{-}+6\Sigma_{\times}A
\end{align}

\begin{flushleft}
\begin{align}
A^{'} & =(q+2\Sigma_{+})A
\end{align}
Abbreviations:
\par\end{flushleft}

\begin{equation}
q=2\Sigma^{2}+\frac{1}{2}(3\gamma-2)\Omega
\end{equation}

\begin{equation}
\Sigma^{2}=\Sigma_{+}^{2}+\Sigma_{-}^{2}+\Sigma_{2}^{2}+\Sigma_{\times}^{2}
\end{equation}
Constraints:

\begin{equation}
\label{con1}
\Omega=1-\Sigma^{2}-N_{-}^{2}-4A^{2}
\end{equation}

\begin{equation}
\label{con2}
g=(\Sigma_{+}+\sqrt{3}\Sigma_{-})A-\Sigma_{\times}N_{-}=0
\end{equation}
\\
Auxilliary Equations:\\
\begin{equation}
\Omega^{'}=[2q-(3\gamma-2)]\Omega
\end{equation}

\begin{equation}
g^{'}=2(q+\Sigma_{+}-1)g
\end{equation}
\\
Note that the auxilliary equations (21) and (22) follow from $(11)-(16)$
and show the invariance of $\Omega=0$ and $g=0$, where $\Omega=0$
results in the vacuum equations. 

We define the phase space for vacuum Bianchi $VI_{-1/9}^{^{*}}$ again by requiring that the constraints (\ref{con1}) and (\ref{con2}) are satisfied. This time, we have six variables and two constraints, yielding again a 4-dimensional state space for the Vacuum models, as in Bianchi IX before.

\begin{defn} (Phase Space for Vacuum Bianchi $VI_{-1/9}^{^{*}}$)
\[
\mathcal{B}=\{\Sigma_{+},\Sigma_{-},\Sigma_{\times},\Sigma_{2},N_{-},A) \; | \; 0=1-\Sigma^{2}-N_{-}^{2}-4A^{2} \; \textup{and} \; g=0\}
\]
\end{defn} 

The equations have been analysed in \cite{hew03}, we give here only a very brief overview about the similarities and differences of Bianchi $VI_{-\frac{1}{9}}^{*}$ compared to Bianchi IX that are relevant for our own research (see chapter \ref{bvi-chapter}).

When we look at the equations, we observe there is also a Kasner circle of fixed points: $\mathcal{K}= \{\Sigma_{\times}=\Sigma_{2}=N_{-}=A=0\}$, leading again to $\Sigma_+^2 + \Sigma_-^2=1$. Similar to Bianchi IX, we can define caps of heteroclinic orbits connecting points of $\mathcal{K}$, but Bianchi $VI_{-\frac{1}{9}}^{*}$ is less symmetric than Bianchi IX. The transitions can be illustrated as follows:

\bigskip
\bigskip

\begin{minipage}[t]{0.3\textwidth}
\includegraphics[width=45mm]{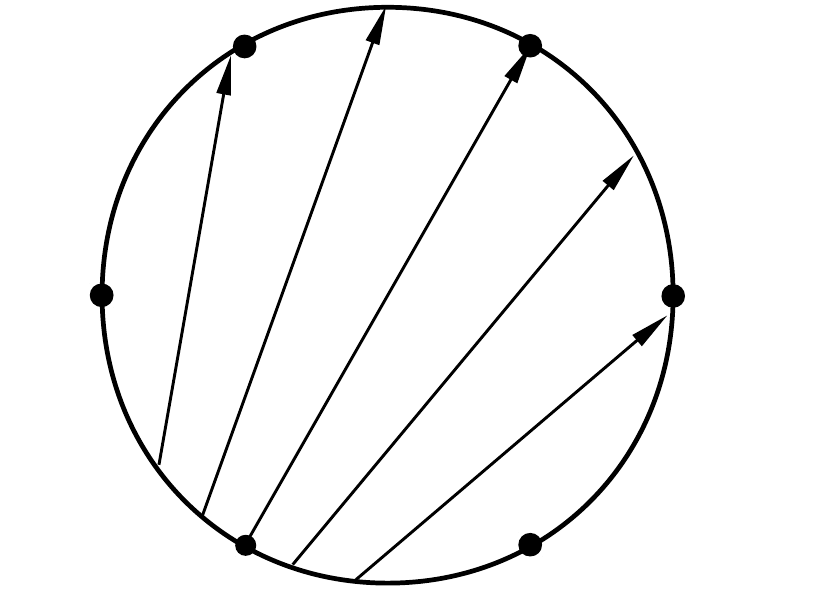}
\center \hspace{0.5cm} $N_-$ 
\end{minipage}
\begin{minipage}[t]{0.3\textwidth}
\includegraphics[width=45mm]{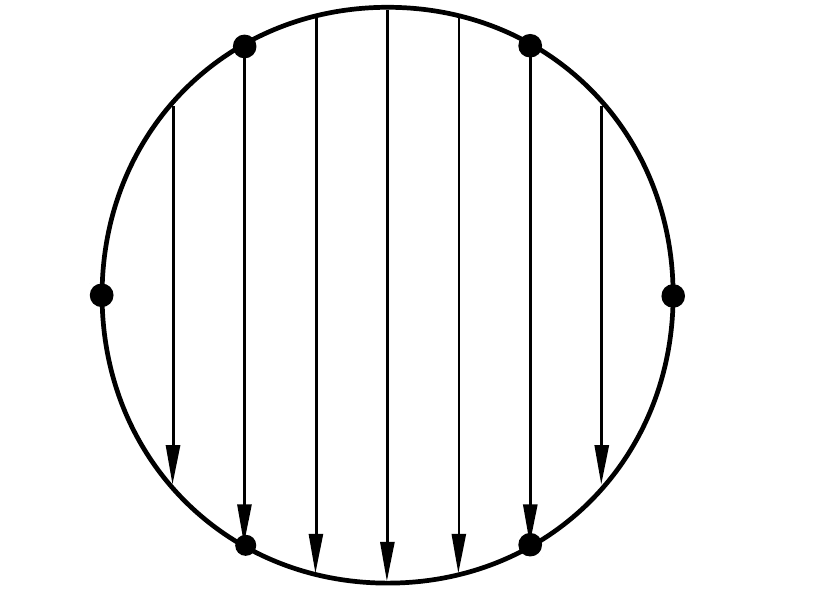}
\center \hspace{0.5cm} $\Sigma_\times$ 
\end{minipage}
\begin{minipage}[t]{0.3\textwidth}
\includegraphics[width=45mm]{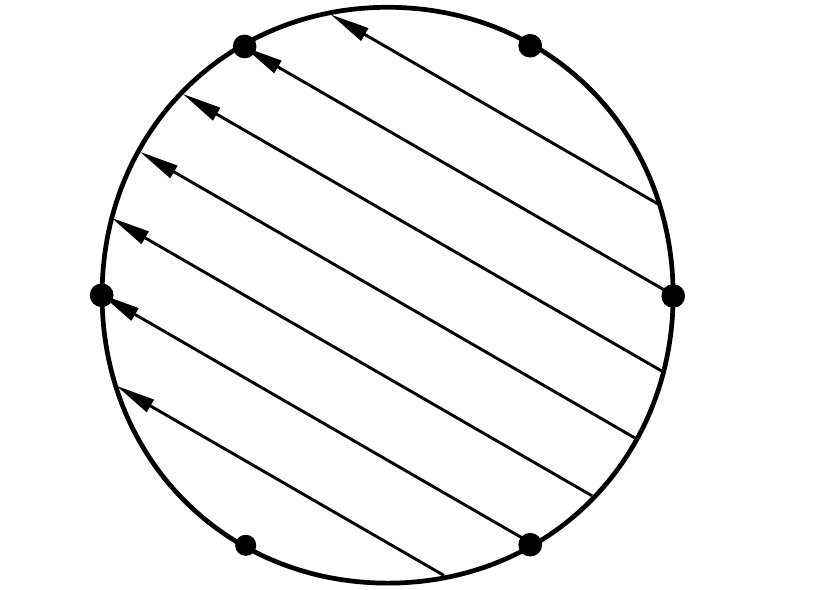}
\center \hspace{0.5cm} $\Sigma_2$ 
\end{minipage}

\bigskip
\bigskip

We have the following caps of heteroclinic orbits: 

\begin{itemize}
\item $C_{N_-}=\{ \Sigma_{\times}=\Sigma_{2}=A=0 \}$ , which represent transitions in the variable $N_{-}$, i.e. curvature transitions. This means the Kasner parameter $u$ changes according to the Kasner map as explained in section \ref{kasner-map} for Bianchi IX.
\item $C_{\Sigma_\times} = \{ \Sigma_{2}=N_{-}=A=0 \}$, which represent transitions in the variable $\Sigma_{x}$, i.e. a frame transition. This means the Kasner parameter $u$ is not changed by the transition, it rather connects two points of $\mathcal{K}$ that are in the same equivalence class with the same $u$ in different sectors
\item $C_{\Sigma_2} = \{ \Sigma_{\times}=N_{-}=A=0 \}$ which represent transitions in the variable in the variable $\Sigma_{2}$, also a frame transition.
\end{itemize}

In addition to the traditional curvature transition in the variable $N_-$ similar to those that also appear in Bianchi class A, we also observe frame transitions of two types, for the variables $\Sigma_{\times}$ and $\Sigma_2$. The fact that the frame transition do not change $u$ can be seen from the fact that the projections of the heteroclinic orbits on the $(\Sigma_+,\Sigma_-)$-plane are parallel lines, and the (inverse) distance to the Taub points (which is one way to interpret $u$) stays the same after the transition. 

\bigskip
\medskip

\begin{minipage}{0.7\textwidth}
\centering
\setlength{\unitlength}{0.277\textwidth}
\begin{picture}(3.6,3.6)(-1.3,-1.8)
\put(-1.3 , 1.8 ){\makebox(0,0)[tl]{%
  \includegraphics[width=3.6\unitlength]{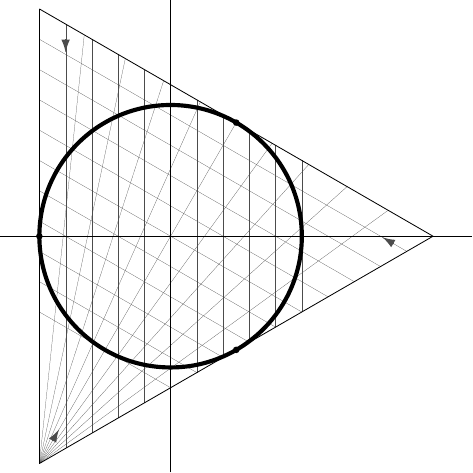}}}
\put( 2.3 ,-0.03){\makebox(0,0)[rt]{$\Sigma_+$}}
\put( 0.03, 1.8){\makebox(0,0)[lt]{$\Sigma_-$}}

\put(-1.03,-0.03){\makebox(0,0)[rt]{$T_1$}}
\put( 0.53, 0.88){\makebox(0,0)[lb]{$T_2$}}
\put( 0.53,-0.88){\makebox(0,0)[lt]{$T_3$}}

\end{picture}
\end{minipage}

\bigskip
\medskip

Similar to what we did in Bianchi IX, we will also find expressions for the eigenvalues of the transition-variables in terms of the Kasner parameter $u$. This will be done in chapter \ref{bvi-chapter}.

\newpage

\section{Existing Results in Bianchi $IX$ and Difficulties in Bianchi B}
\label{DiffBianchiB}

Historically, the chaotic oscillations in Bianchi models have first been observed in the Bianchi IX model, see \cite{lk63,bkl70,bkl82}. It is also known under the name "Mixmaster", a term coined by Misner \cite{mis69a}. The first rigorous theorem on the ancient dynamics in Bianchi IX was proved by Ringström \cite{rin01}, based on earlier results by Rendall \cite{ren97}. For a recent survey on the "Facts and Beliefs" concerning the Mixmaster model, see \cite{heiugg09b}.


One important research question is relating properties of the Kasner map (which is known to be chaotic, see e.g. the chapter 11 
 in the book \cite{waiell97}) to the real dynamics in Bianchi models. This can be seen as making the BKL-conjecture rigorous for spatially homogenous spacetimes, i.e. by proving the "oscillatory" (and possibly also the "vacuum") part in a case where the model is already "local".

In Bianchi IX, there exist such rigorous convergence results by Liebscher et al \cite{lieetal10}, B\'eguin \cite{beg10}, and Reiterer/Trubowitz \cite{reitru10}. For Bianchi  $VI_{0}$ with a magnetic field as matter, there have been recent results by Liebscher, Rendall and Tschapa \cite{lieetal12}. 


Until today, there exist no rigorous convergence results for Bianchi $VI_{-1/9}^{^{*}}$, which is of class B (see chapter \ref{gen-rel}, section \ref{symm-sh-models}).
The reason for this is that Bianchi $VI_{-1/9}^{^{*}}$ is more difficult than Bianchi class A. One example is that in Bianchi IX, the normal hyperbolicity of the Kasner circle fails only at the three Taub points, while in Bianchi $VI_{-1/9}^{^{*}}$ this is true for all of the six points that mark the borders of the sectors of the Kasner circle defined in section \ref{bianchi-defs}. Also there are non-unique heteroclinic chains because of multiple unstable eigenvalues for some sectors of the Kasner circle, marked in red in the picture below. We will discuss this matter further in chapter \ref{bvi-chapter}. 

\bigskip
\bigskip
\begin{minipage}{0.4\textwidth}
\centering
\setlength{\unitlength}{0.277\textwidth}
\begin{picture}(3.6,3.6)(-1.3,-1.8)
\put(-1.3 , 1.8 ){\makebox(0,0)[tl]{%
  \includegraphics[width=3.6\unitlength]{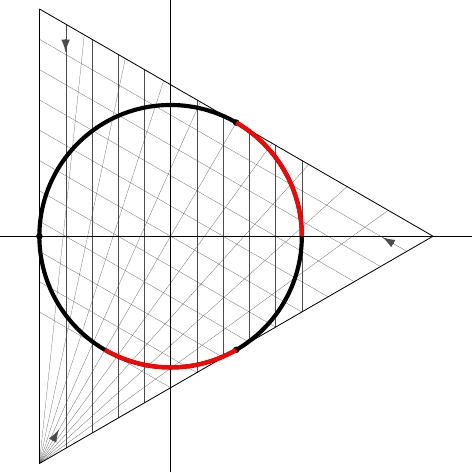}}}
\put( 2.3 ,-0.03){\makebox(0,0)[rt]{$\Sigma_+$}}
\put( 0.03, 1.8){\makebox(0,0)[lt]{$\Sigma_-$}}

\put(-1.03,-0.03){\makebox(0,0)[rt]{$T_1$}}
\put( 0.53, 0.88){\makebox(0,0)[lb]{$T_2$}}
\put( 0.53,-0.88){\makebox(0,0)[lt]{$T_3$}}

\end{picture}
\end{minipage}



\newpage

\newpage

\chapter{Combined Linear Local Passage at the 18-cycle in Bianchi $VI_{-1/9}^{^{*}}$}
\label{bvi-chapter}

In this section, we show how some of the techniques developed in the chapter before can be applied to Bianchi $VI_{-1/9}^{^{*}}$, where no rigorous convergence result exists to date. 
We construct an example for a periodic heteroclinic
chain in Bianchi $VI_{-1/9}^{^{*}}$ that allows Takens Linearization
at all Base points. It will turn out to be a "18-cycle'', i.e.
involving a heteroclinic chain of 18 different base points at the
Kasner circle. We then show that the combined linear local passage
at the 18-cycle is a contraction. This qualifies
it as a possible candidate for proving a rigorous convergence
theorem in Bianchi $VI_{-1/9}^{^{*}}$.

The situation in Bianchi $VI_{-1/9}^{^{*}}$ is more involved than in Bianchi IX, as there are sectors of the Kasner circle with more than one unstable eigenvalue (towards the big bang, the time direction we are interested in). This means that even by starting with the same Kasner-parameter, one can have different realizations in terms of  heteroclinic chains. 

If we label the six sectors  of the Kasner circle counter-clockwise as it is often done in Bianchi models (starting in the positive quadrant), this is true e.g. for sector 5, which is part of the 18-cycle. In the graphics below, the sectors with multiple unstable eigenvalues are marked in red color. The different families of heteroclinic orbits in Bianchi $VI_{-1/9}^{^{*}}$ are illustrated by the lines in light gray in the background (see chapter \ref{bianchi-defs}, section \ref{bianchi-vi-defs} for details):

\medskip

\begin{minipage}[t]{0.2\textwidth}
\end{minipage}
\begin{minipage}{0.4\textwidth}
\centering
\setlength{\unitlength}{0.277\textwidth}
\centering
\begin{picture}(3.6,3.6)(-1.3,-1.8)
\put(-1.3 , 1.8 ){\makebox(0,0)[tl]{%
  \includegraphics[width=3.6\unitlength]{pics/BianchiB-intro}}}
\put( 2.3 ,-0.03){\makebox(0,0)[rt]{$\Sigma_+$}}
\put( 0.03, 1.8){\makebox(0,0)[lt]{$\Sigma_-$}}

\put(-1.03,-0.03){\makebox(0,0)[rt]{$T_1$}}
\put( 0.53, 0.88){\makebox(0,0)[lb]{$T_2$}}
\put( 0.53,-0.88){\makebox(0,0)[lt]{$T_3$}}

\put(1.03,-0.03){\makebox(0,0)[lt]{$Q_1$}}
\put( -0.53, 0.88){\makebox(0,0)[rb]{$Q_2$}}
\put( -0.53,-0.88){\makebox(0,0)[rt]{$Q_3$}}

\put( 0,1.03){\makebox(0,0)[lb]{{\bf \Large 2 }}}
\put( 0,-1.03){\makebox(0,0)[lt]{{\bf \Large 5 }}}
\put( 0.88, 0.53){\makebox(0,0)[bl]{{\bf \Large 1 }}}
\put( -0.88,0.53){\makebox(0,0)[rb]{{\bf \Large 3 }}}
\put( -0.88, -0.53){\makebox(0,0)[tr]{{\bf \Large 4}}}
\put( 0.88,-0.53){\makebox(0,0)[lt]{{\bf \Large 6 }}}
\end{picture}
\end{minipage}
\begin{minipage}[t]{0.2\textwidth}
\end{minipage}

We will illustrate this ambiguity that can arise from the same continued fraction development $u={[}3,5,3,5,...{]}$ by showing two possible "18-cycles" (details are discussed in section \ref{PassagesBianchiB}):

\begin{itemize}
\item the ``classic'' 18-cycle, with sector sequence 54343-425-43434343-425,
where we go ``left'' to sector 4 both times in sector 5
\item the ``advanced'' 18-cycle, with sequence 5634342-543434343425,
where we go first ``right'' to sector 6, and the second time ``left''
to sector 4
\end{itemize}

The graphics below shows the "advanced 18 cycle" starting in sector 5. The bounces around the Taub-point can be seen clearly - the blue color shows the part with 3 bounces, while the purple
color exhibits 5 bounces. The label "advanced 18-cylce" means that we first follow the blue part of the heteroclinic chain and then the purple one, reflecting our choice of  u={[}3,5,3,5,...{]}.

\medskip
\begin{minipage}{0.6\textwidth}
\centering
\setlength{\unitlength}{0.277\textwidth}
\begin{picture}(3.6,3.6)(-1.3,-1.8)
\put(-1.3 , 1.8 ){\makebox(0,0)[tl]{%
  \includegraphics[width=3.6\unitlength]{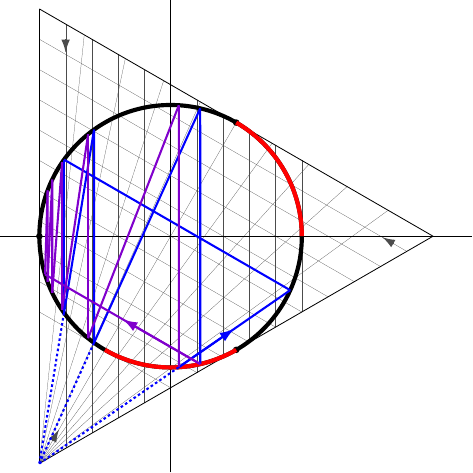}}}
\put( 2.3 ,-0.03){\makebox(0,0)[rt]{$\Sigma_+$}}
\put( 0.03, 1.8){\makebox(0,0)[lt]{$\Sigma_-$}}

\put(-1.03,-0.03){\makebox(0,0)[rt]{$T_1$}}
\put( 0.53, 0.88){\makebox(0,0)[lb]{$T_2$}}
\put( 0.53,-0.88){\makebox(0,0)[lt]{$T_3$}}

\end{picture}
\end{minipage}
\medskip



Now we are in a position to state the two main results of this chapter:

\begin{thm}
\label{thm-18-takens}
Both the classic and the advanced 18-cycle allow Takens-Linearization
at all base points. 
\end{thm}

\begin{lem}
\label{thm-18cllp}
The Combined Linear Local Passage at the classic 18-cycle is a contraction.
\end{lem}

The proof of these results is the aim of this chapter, and we will now comment a little on how the rest of  chapter is organized. 

In order to illustrate our approach, we start by considering the 3-cycle in Bianchi $VI_{-1/9}^{^{*}}$, as it is easier to handle and the method we develop applies also to longer cycles.

First we derive the formulas for the relevant eigenvalues of the linearized vector field at the base-points of the 3-cycle  in Bianchi $VI_{-1/9}^{^{*}}$ in terms of the Kasner Parameter u. We then consider the Combined Linear Local Passages and Takens Linearization at the 3-cycle. We are able to show that the Sternberg Non-Resonance
Conditions are not satisfied for the 3-cycle, but they are satisfied for the 18-cycles. 

This means that both the classic and the advanced 18-cycles are periodic heteroclinic chains in Bianchi $VI_{-1/9}^{^{*}}$ that allow Takens Linearization at all of their Base points. In addition, we give evidence that the Combined Linear Local Passage at the classic 18-cycle is a contraction. This qualifies it as a candidate for proving a rigorous convergence theorem in Bianchi $VI_{-1/9}^{^{*}}$. 

In order to progress further and turn Lemma \ref{thm-18cllp} into a rigorous convergence theorem, a better understanding of the global passages in Bianchi $VI_{-1/9}^{^{*}}$ is necessary. We will comment on possible ideas how to do this in the section "Conclusion and Outlook".

As the formulas for the eigenvalues at points on the Kasner circle in Bianchi $VI_{-1/9}^{^{*}}$ depend on the sector, many cases have to be checked. That's why we use Mathematica in order to do the necessary computations for the 18-cycles (see Appendix \ref{Appendix18cycle}).

\newpage

\section{Eigenvalues in Terms of the Kasner Parameter u}

\subsection{General Formulas for Points on the Kasner Circle}

At first, we recall that in Bianchi IX, the formula for the eigenvalues of the linearized vector field at points of the Kasner circle have an easy expression in terms of the Kasner parameter u, compare  section \ref{bianchi-defs}, equation (\ref{eq:EVinBIXinU}). In Bianchi IX, the situation in different sectors of the Kasner-circle only differs by a permutation of those 3 formulas for the eigenvalues (see \cite{heiugg09b}, p.8), which does not matter for the question of resonances. 

But in Bianchi $VI_{-1/9}^{^{*}}$, the situation is more complicated. Here, the formulas for the eigenvalues at points on the Kasner circle do depend
on the sector, so in order to check for resonances, a lot of cases have to be considered. Each sector corresponds to a permutation of the Kasner expontents $p_{i}$,
where sector (321) means e.g. that $p_{3}<p_{2}<p_{1}$which fixes
the formula for each of them. For this is it important to note that
$u\in[1,\infty]$. As an example, consider the sector 5 or (312), compare section \ref{bianchi-defs}): 
\[
p_{3}=\frac{-u}{1+u+u^{2}}
\]

\[
p_{1}=\frac{(u+1)}{1+u+u^{2}}
\]

\[
p_{2}=\frac{u(u+1)}{1+u+u^{2}}
\]
The general formuals for the relevant eigenvalues (i.e. $\lambda_{\times}$, $\lambda_{2}$, $\lambda_{-}$ corresponding to the variables involved in the heteroclinic chain: $\Sigma_{\times}$, $\Sigma_{2}$, $N{-}$, see section \ref{bianchi-vi-defs}) in terms of $\Sigma_{+},\Sigma_{-}$resp.
the $p_{i}$ are (see \cite{heiugg09b}, p.7):

\[
\Sigma_{+}=-\frac{3}{2}p_{1}+\frac{1}{2}
\]

\[
\Sigma_{-}=-\frac{\sqrt{3}}{2}(p_{1}+2p_{2}-1)
\]

\[
\lambda_{\times}=-2\sqrt{3}\Sigma_{-}
\]

\[
\lambda_{2}=-3\Sigma_{+}+\sqrt{3}\Sigma_{-}
\]

\[
\lambda_{-}=2+2\Sigma_{+}+2\sqrt{3}\Sigma_{-}
\]

\newpage

\subsection{Eigenvalues at the 3-Cycle}
\label{ev-3cyle}

In the following, we present the formulas for the sectors that are involved
in the 3-cycle with $u=\mbox{golden mean}=\frac{1+\sqrt{5}}{2}$ and
the sector-sequence ``5-1-2-5''. For the other sectors, similar
formulas can be derived analogously (see Appendix \ref{CheckTakens18cycle}).


\bigskip

\paragraph*{Base Point $B_{1}$ in sector 5, i.e. (312) }

\bigskip

\[
\lambda_{2}=\frac{3-3u^{2}}{1+u+u^{2}}
\]

\[
\lambda_{\times}=\frac{6u+3u^{2}}{1+u+u^{2}}
\]

\[
\lambda_{-}=\frac{-6u}{1+u+u^{2}}
\]

\bigskip

\paragraph*{Base Point $B_{2}$ in sector 1, i.e. (123) }

\bigskip

\[
\lambda_{2}=\frac{-3-6u}{1+u+u^{2}}
\]

\[
\lambda_{\times}=\frac{3-3u^{2}}{1+u+u^{2}}
\]

\[
\lambda_{-}=\frac{6u+6u^{2}}{1+u+u^{2}}
\]

\bigskip

\paragraph*{Base Point $B_{3}$ in sector 2, i.e. (213) }

\bigskip

\[
\lambda_{2}=\frac{3+6u}{1+u+u^{2}}
\]

\[
\lambda_{\times}=\frac{-6u-3u^{2}}{1+u+u^{2}}
\]

\[
\lambda_{-}=\frac{6u+6u^{2}}{1+u+u^{2}}
\]

\newpage

\section{The 3-Cycle in Bianchi $VI_{-1/9}^{^{*}}$}

\subsection{(Non-)Resonance and Takens Linearization }

In this section, we shortly check the Sternberg-Non-Resonance-Conditions
(SNC) for the 3-cycle in Bianchi $VI_{-1/9}^{^{*}}$. The procedure
necessary to do this is described in detail in section \ref{resonances-bix}.
The parameter-value at the 3-cylce is $u=g=\frac{1+\sqrt{5}}{2}$,
which satisfies 

\[
1+\frac{1}{u}=u\implies1+u-u^{2}=0
\]

The equation for checking the (SNC) thus reads:

\[
M*\left(\begin{array}{c}
k_{1}\\
k_{2}\\
k_{3}
\end{array}\right)=z*\left(\begin{array}{c}
1\\
1\\
-1
\end{array}\right)
\]

According to the formulas above, we observe that

\[
M_{B_{1}}=\left(\begin{array}{ccc}
3 & 0 & 0\\
0 & 6 & -6\\
\text{-3} & 3 & 0
\end{array}\right),M_{B_{2}}=\left(\begin{array}{ccc}
-3 & 3 & 0\\
-6 & 0 & 6\\
\text{0} & -3 & 6
\end{array}\right),M_{B_{3}}=\left(\begin{array}{ccc}
3 & 0 & 0\\
6 & -6 & 6\\
\text{0} & -3 & 6
\end{array}\right)
\]
, which are all invertible, and give the following results (with $z=6$
for the earliest possible resonance):

\[
k_{B_{1}}=\left(\begin{array}{c}
2\\
0\\
-1
\end{array}\right),k_{B_{2}}=\left(\begin{array}{c}
-2\\
0\\
-1
\end{array}\right),k_{B_{3}}=\left(\begin{array}{c}
2\\
0\\
-1
\end{array}\right)
\]

This means that (SNC) does hold at the Base points of the 3-cycle
up to order 3, which is not enough to allow Takens-Linearization,
as the the order necessary in the Takens-Theorem, named $\alpha$(1),
is bigger than 10 in all of the sectors involved (see Appendix \ref{CheckTakens3cycle}).
That's why Takens Linearization Theorem may not be employed at the
base points of the 3-cycle, and we have to look for a different (longer)
cycle. Nevertheless, we will stick to the 3-cycle in the following
chapter in order to illustrate our method for calculating the Combined
Linear Local Passage, as the method applies to longer cycles as well.

\newpage

\subsection{Combined Linear Local Passages}
\label{3cycle-cllp}

In this section, we use the linearized vectorfield at the base points
of the 3-cycle to explicitly compute the local passages as was done for Bianchi IX before (see \ref{compute-LP-start}-\ref{compute-LP-end}). Be aware that Takens-Linearization is not allowed at the 3-cycle in Bianchi $VI_{-1/9}^{^{*}}$, so this is only a formal calculation in this case to illustrate what we mean by "Combined Local Linear Passage". Later, for the 18-cycles, the calculation will be justified as Takens-Linearization is possible there.

The ratio of the relevant Eigenvalues at the Base-points $B_{1}$,$B_{2}$,$B_{3}$ of the 3-cycle
(corresponding to the value $u=\frac{1+\sqrt{5}}{2}$ in the 3 sectors)
is given by the following:\\

$\mbox{Near }\ensuremath{B_{1}=(-\frac{1}{4},-\frac{\sqrt{15}}{4}):}$

\[
r_{1}=-\frac{\lambda_{\times}}{\lambda_{n}}=\frac{u+2}{2}=1.8090>0
\]

\[
r_{2}=-\frac{\lambda_{2}}{\lambda_{n}}=-\frac{u^{2}-1}{2u}=-0.5000<0
\]

\bigskip

$\mbox{Near }\ensuremath{B_{2}=(\frac{1+3\sqrt{5}}{8},\frac{\sqrt{15}-\sqrt{3}}{8}):}$

\[
r_{3}=-\frac{\lambda_{n}}{\lambda_{2}}=\frac{2u(u+1)}{2u+1}=2>0
\]

\[
r_{4}=-\frac{\lambda_{u}}{\lambda_{uu}}=-\frac{\lambda_{\times}}{\lambda_{2}}=-\frac{u^{2}-1}{2u+1}=-0.3820<0
\]

\bigskip

$\mbox{Near }\ensuremath{B_{3}=(-\frac{1}{4},\frac{\sqrt{15}}{4}):}$

\[
r_{5}=-\frac{\lambda_{2}}{\lambda_{\times}}=\frac{2u+1}{u(u+2)}=0.7236>0
\]

\[
r_{6}=-\frac{\lambda_{n}}{\lambda_{\times}}=\frac{2(u+1)}{u+2}=1.4472>0
\]

\bigskip

We start our combined linear local passage at the In-Section of the local passage at $B_{1}$ in sector 5. Thus $\Sigma_{\times}$ is the incoming variable and we define $a:=\Sigma_{2}^{in}$, $b:=N_{-}^{in}$ (the two remaining relevant variables in the section).

Then we do the calculations for the 3 local passages involved at the 3-cycle near the $B_{1}$,$B_{2}$ and $B_{3}$ 
as in Bianchi IX before (see \ref{compute-LP-start}-\ref{compute-LP-end}). We arrive at the following formulas for the combined linear local passage near the 3-cycle:

\[
\tilde{a}=([b^{r_{2}}a]^{r_{4}}b^{r_{1}})^{r_{5}}
\]

\[
\tilde{b}=([b^{r_{2}}a]^{r_{4}}b^{r_{1}})^{r_{6}}(b^{r_{2}}a)^{r_{3}}
\]

Taking Logarithms on both sides yields the following:\\

\[
\left(\begin{array}{c}
\log\tilde{a}\\
\log\tilde{b}
\end{array}\right)=M_{3-cycle}*\left(\begin{array}{c}
\log a\\
\log b
\end{array}\right)
\]
\\

with the following matrix $M_{3-cycle}$:

\[
M_{3-cycle}=\left(\begin{array}{cc}
r_{5}r_{4} & r_{5}r_{4}r_{2}+r_{5}r_{1}\\
r_{6}r_{4}+r_{3} & r_{6}r_{4}r_{2}+r_{6}r_{1}+r_{3}r_{2}
\end{array}\right)
\]

\bigskip

As mentioned before, Takens-Linearization is not allowed at the 3-cycle in Bianchi $VI_{-1/9}^{^{*}}$.
However, for the 18-cycles, an analogous calculation will be justified as Takens-Linearization is possible there. We will discuss the properties $M_{18-cycle}$ in order to prove Theorem \ref{thm-18cllp} in section \ref{18cycle-cllp}.


\newpage

\section{The 18 Cycles in Bianchi $VI_{-1/9}^{^{*}}$}

\subsection{Possible Passages in Bianchi $VI_{-1/9}^{^{*}}$} \label{PassagesBianchiB}

When we look at the different transitions possible in Bianchi $VI_{-1/9}^{^{*}}$ (see chapter \ref{bianchi-defs}, section \ref{bianchi-vi-defs}), 
we are able to understand which sequence of sectors can arise when solutions converge to their corresponding heteroclinic chains. For the classification below, we have chosen to put our section always before the next ``Curvature Transition'', i.e. when we leave from sector 4 or 5,
that's why all the passages start and end in one of these sectors:

\begin{itemize}
\item Passage A: Sectors 4-3-4
\item Passage B1: Sectors 4-2-5
\item Passage B2: Sectors 4-2-5-4
\item Passage C1: Sectors 5-1-2-5
\item Passage C2: Sectors 5-1-2-5-4
\item Passage D: Sectors 5-6-3-4
\item Passage E: Sectors 5-1-6-3-4
\end{itemize}
For the 18-cycles discussed below, only a few of the passages will occur, namely A, B1, B2 and D. 



\subsection{The Classic 18-Cycle}

We now start with u={[}3,5,3,5,...{]} in Sector 4 and prescribe the follwing dynamics:\\

\begin{tabular}{|c|c|c|c|c|c|c|c|}
\hline 
u= & {[}3,5,...{]} & {[}2,5,...{]} & {[}2,5,...{]} & {[}1,5,...{]} & {[}1,5,...{]} & {[}5,3,...{]} & {[}5,3,...{]}\tabularnewline
\hline 
sector & 4 & 3 & 4 & 3 & 4 & 2 & 5\tabularnewline
\hline 
\end{tabular}\\

\begin{tabular}{|c|c|c|c|c|c|c|c|}
\hline 
u= & {[}5,3,...{]} & {[}4,3,...{]} & {[}4,3,...{]} & {[}3,3,...{]} & {[}3,3,...{]} & {[}2,3,...{]} & {[}2,3,...{]}\tabularnewline
\hline 
sector & 4 & 3 & 4 & 3 & 4 & 3 & 4\tabularnewline
\hline 
\end{tabular}\\

\begin{tabular}{|c|c|c|c|c|c|c|c|}
\hline 
u= & {[}1,3,...{]} & {[}1,3,...{]} & {[}3,5,...{]} & {[}3,5,...{]} & {[}3,5,...{]} & ... & ...\tabularnewline
\hline 
sector & 3 & 4 & 2 & 5 & 4 & ... & ...\tabularnewline
\hline 
\end{tabular}\\
\\
This means we have the following sequence of Passages: A A B2 A A
A A B2, and this pattern continues arbitrarily often. This involves
18 global passages, that's why I call it an "18-cycle''. 

Observe that both 18-cycles mentioned in the introduction to this chapter started in sector 5.
This was done for illustrative purposes as there is an ambiguity how to continue in this sector.
From now on, we refer to the ``classic 18-cycle'' with the sequence of sectors as illustrated in the table above, which means we have started in sector 4.

Note that we could also derive a different sequence of passages for the same
u, as we have a choice in sector 5 either to go to sector 4 via a
frame transition (as done above) or to go via curvature transition
to sector 6, as done at the first transition for the advanced 18-cycle in the next section.

\subsection{The Advanced 18-Cycle}

Now we start with u={[}3,5,3,5,...{]} in Sector 5 and prescribe the follwing
dynamics:\\

\begin{tabular}{|c|c|c|c|c|c|c|c|}
\hline 
u= & {[}3,5,...{]} & {[}2,5,...{]} & {[}2,5,...{]} & {[}2,5,...{]} & {[}1,5,...{]} & {[}1,5,...{]} & {[}5,3,...{]}\tabularnewline
\hline 
sector & 5 & 6 & 3 & 4 & 3 & 4 & 2\tabularnewline
\hline 
\end{tabular}\\

\begin{tabular}{|c|c|c|c|c|c|c|c|}
\hline 
u= & {[}5,3,...{]} & {[}5,3,...{]} & {[}4,3,...{]} & {[}4,3,...{]} & {[}3,3,...{]} & {[}3,3,...{]} & {[}2,3,...{]}\tabularnewline
\hline 
sector & 5 & 4 & 3 & 4 & 3 & 4 & 3\tabularnewline
\hline 
\end{tabular}\\

\begin{tabular}{|c|c|c|c|c|c|c|c|}
\hline 
u= & {[}2,3,...{]} & {[}1,3,...{]} & {[}1,3,...{]} & {[}3,5,...{]} & {[}3,5,...{]} & ... & ...\tabularnewline
\hline 
sector & 4 & 3 & 4 & 2 & 5 & ... & ...\tabularnewline
\hline 
\end{tabular}\\
\\
This means we have the following sequence of passages: D A B2 A A
A A B1, which defines the "advanced 18-cycle".

\subsection{(Non-)Resonance and Takens Linearization at the 18-Cycle}

We now show that both of the 18 cycles are infinite periodic
heteroclinic chains in Bianchi $VI_{-1/9}^{^{*}}$ that allows Takens
Linearization at all of its base points. The Mathematica-output in Appendix \ref{CheckTakens18cycle}
shows that for $u=[3,5,3,5,...]$ the Sternberg-Non-Resonance-Conditions
(SNC) are satisfied, because the $\alpha(1)$ that is necessary for a $C^{1}$-linearization at
each point is always smaller than the sum of the absolute value of the
coefficients in the vector $k=\{k_{1},k_{2},k_{3}\}.$ That's why
we can employ the Takens Linearization Theorem for both of the 18-cycles
mentioned above, and Theorem \ref{thm-18-takens} is proven.

\subsection{Combined Linear Local Passage at the 18-Cycle}
\label{18cycle-cllp}

Our Results from Mathematica (see Appendix \ref{CLLP-18cycle}) indicate that we get the following matrix for the
combined linear local passage of the classic 18-cycle when we apply the same algorithm that we outlined for the 3-cycle in section \ref{3cycle-cllp}:

\[
M_{18-cycle}=\left(\begin{array}{cc}
267.54 & 110.78\\
595.16 & 247.51
\end{array}\right)=:\left(\begin{array}{cc}
v_{1} & v_{2}\\
v_{3} & v_{4}
\end{array}\right)
\]

\[
\mu_{1}=514.49\mbox{ with Eigenvector }v_{1}=(0.45,1)
\]

\[
\mu_{2}=0.5578\mbox{ with Eigenvector }v_{2}=(-0.41,1)
\]

\bigskip

This implies that if we start with small positive a,b (i.e. log a,
log b <\textcompwordmark{}< 0), the combined linear local passage
will bring us closer to the origin - this is what we mean by the term "contraction" in Lemma \ref{thm-18cllp}.

\newpage

\section{Numerical Simulation}

The picture below shows a numerical simulation (with Matlab) of a periodic heteroclinic chain in Bianchi $VI_{-1/9}^{^{*}}$, here a 13-cycle. According to the terminology developed above, it would be named the "advanced" 13 cycle, as both directions are taken from sector 5.  

Although a much more detailed numerical analysis is necessary, our simulation shows that at least there are cases where both directions are taken from sector 5. Thus this possibility seems to really occur in the equations, at least numerically there are Bianchi $VI_{-1/9}^{^{*}}$-solutions following the "advanced" 13 cycle towards the big bang.

Of course much more effort is needed in order to set up the numerics in an appropriate way instead of just using a built-in Matlab ODE solver\footnote{for producing the picture above, we have used the "ode113" solver, which is a variable order Adams-Bashforth-Moulton PECE solver (according to the Matlab documentation \cite{matlab-ode}). We thank Woei Chet Lim for providing us with some Matlab code that we used in order to carry out our numerical simulations.}. One idea could be to use the explicit linear flow near the equilibria of the Kasner circle, where most of the time is spent, in order to achieve a higher precision.

\bigskip

\medskip
\begin{center}
\includegraphics[scale=0.75]{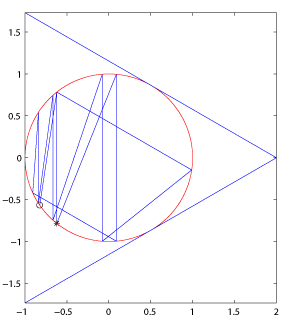}
\end{center}
\medskip

\newpage

\chapter*{Conclusion and Outlook}
\label{conl-outl}

We have constructed, for the fist time, the 18-cycle as a relatively simple example for a periodic heteroclinic chain in Bianchi $VI_{-1/9}^{^{*}}$ that allows Takens Linearization at all Base points.  In addition, we were able to show that the Combined Linear Local Passage at the classic 18-cycle
is a contraction. This could be seen as a first step for proving a rigorous convergence theorem in Bianchi $VI_{-1/9}^{^{*}}$.

In order to progress further, a better understanding of the global passages in Bianchi $VI_{-1/9}^{^{*}}$ is necessary. This is not an easy task, as there are less invariant subspaces than in Bianchi IX that restrict the signs of the heteroclinic orbits, so much more complicated transitions are possible. A first step could be to check in detail which sequence of signs for the different transitions occurs numerically, leading to a classification of possible cases. Speculatively, one could think about the possibility to prove a theorem that for heteroclinic chains with periodic continued fraction developments and with a clearly defined seqence of signs for the transitions there exist solutions of the Bianchi $VI_{-1/9}^{^{*}}$ equations that show this behaviour. But this matter requires further investigation.

Until now, the application of Dynamical Systems Techniques to spatially homogeneous cosmological models yielding Ordinary Differential Equations has been discussed. However, there has also been the attempt to apply such techniques to inhomogeneous cosmologies, yielding Partial Differential Equations. The reason is that in a way the main point\footnote{let us again quote the survey paper by Uggla on this issue (\cite{ugg13a},p.2): "However, arguably the most central, and controversial, assumption of BKL is their `locality' conjecture. According to BKL, asymptotic dynamics toward a generic spacelike singularity in inhomogeneous cosmologies is `local,' in the sense that each spatial point
is assumed to evolve toward the singularity individually and independently of its
neighbors as a spatially homogeneous model"} of the BKL-picture is the question of "locality", asking if the "complicated" Einstein Equations that are PDEs\footnote{for General Relativity from the viewpoint of Partial Differential Equations see \cite{ren08,rin09}} can be approximated by "simpler/less complicated" ODEs towards the big bang. Until today, mostly numerical and heuristic results exist in inhomogeneous cosmologies, but very few rigorous mathematical theorems. 

In the paper "The past attractor in inhomogeneous cosmologies" (\cite{uggetal03}), it is outlined how 
it could be achieved to make the "local" part of the BKL-picture more rigorous, compare also \cite{heietal09}. After some results have been proven in the oscillatory spatially homogeneous setting of
Bianchi IX and in an inhomogeneous, but non-oscillatory setting of Gowdy-spacetimes (see e.g. \cite{ren08}), the logical next step seems to be to consider inhomogeneous oscillatory cosmological models. Arguably the simplest case is given by the $G_2$-cosmologies, that's why it has received rising attention in recent years (\cite {elsetal02,lim04,limetal09,buchner-evolution}). However, there is not a single rigorous convergence theorem comparable to the results that could be achieved in spatially homogeneous models.

A particular complication in inhomogeneous models is the occurrence of spikes, i.e. the formation of
spatial structure. Numerical experiments support the conjecture that spikes form the non-local part of
the generalized Mixmaster attractor (\cite{limetal09}). Lim has also found explicit spike solutions that are compatible with the usual Bianchi II - transitions, giving rise to a "non-local" version of the mimaster dynamics involving "spike transitions" (\cite{lim08}). Recent progress has been achieved by Heinzle and Uggla, who report more in detail about the role of the spike solutions as building blocks of such an extended non-local mixmaster-dynamics (\cite{heietal12}). In addition, they have done a statistical analysis on the spikes in $G_{2}$-models (\cite{heiugg13}). 

A first step towards achieving rigorous results in inhomogeneous oscillatory models could  be to investigate the process of spike formation in $G_{2}$-models. A good understanding of the underlying spatially homogeneous model (which is Bianchi $VI_{-1/9}^{^{*}}$) is probably necessary for this project, but as Uggla writes in his recent survey: "Unfortunately, there exist no rigorous mathematical results concerning their past asymptotic dynamics", refering to Bianchi $VI_{-1/9}^{^{*}}$ (\cite{ugg13a}, p.11).

\newpage

\appendix


\section{Results on Non-Resonance-Conditions and CLLP for Heteroclinic Cycles in Bianchi $VI_{-1/9}^{^{*}}$} 

\label{Appendix18cycle} 


\subsection{Takens Linearization at the Base Points of the 3-Cycle}    
\label{CheckTakens3cycle}

 For u=[m,a,a,...] and m=1...a with the following PARAMETERS: \\
 a= 1, k= 1 (k=1 means smoothness is  $C^{1}$)

 m= 1

 Sektor 1  alpha=  22  beta= 3  k= {2,0,1}

 Sektor 2  alpha=  14  beta= 5  k= {-2,0,1}

 Sektor 3  alpha=  9  beta= 3  k= {2,4,1}

 Sektor 4  alpha=  15  beta= 6  k= {2,0,1}

 Sektor 5  alpha=  16  beta= 3  k= {-2,0,1}

 Sektor 6  alpha=  17  beta= 8  k= {2,4,1}

\subsection{Takens Linearization at the Base Points of the 18-Cycle}
\label{CheckTakens18cycle} We first give the coefficient matricies for the eigenvalues in the other sectors not part of the 3-cycle (see chapter \ref{bvi-chapter}, section \ref{ev-3cyle}):

\[
M_{S_{3}}=\left(\begin{array}{ccc}
0 & -3 & 6\\
6 & -6 & 6\\
3 & 0 & 0
\end{array}\right),M_{S_{4}}=\left(\begin{array}{ccc}
-3 & 3 & 0\\
0 & 6 & -6\\
3 & 0 & 0
\end{array}\right),M_{S_{6}}=\left(\begin{array}{ccc}
0 & -3 & 6\\
-6 & 0 & 6\\
-3 & 3 & 0
\end{array}\right)
\]
For u={[}m,a,b,a,b,...{]} and m=1...b with the following 
PARAMETERS: a= 3 b= 5 k= 1 (k=1 means smoothness is $C^{1}$)\\

m= 1

Sektor 1 alpha= 34 beta= 3 k= \{-46,-80,-37\}

Sektor 2 alpha= 11 beta= 4 k= \{-34,-80,-37\}

Sektor 3 alpha= 8 beta= 3 k= \{6,-40,-37\}

Sektor 4 alpha= 25 beta= 10 k= \{6,-28,-37\}

Sektor 5 alpha= 25 beta= 3 k= \{-34,-28,-37\}

Sektor 6 alpha= 30 beta= 14 k= \{-46,-40,-37\}

m= 2

Sektor 1 alpha= 15 beta= 3 k= \{-46,-92,-43\}

Sektor 2 alpha= 17 beta= 6 k= \{-46,-92,-43\}

Sektor 3 alpha= 13 beta= 4 k= \{6,-40,-43\}

Sektor 4 alpha= 10 beta= 4 k= \{6,-40,-43\}

Sektor 5 alpha= 11 beta= 3 k= \{-46,-40,-43\}

Sektor 6 alpha= 10 beta= 5 k= \{-46,-40,-43\}

m= 3

Sektor 1 alpha= 16 beta= 3 k= \{-34,-80,-37\}

Sektor 2 alpha= 22 beta= 8 k= \{-46,-80,-37\}

Sektor 3 alpha= 15 beta= 4 k= \{6,-28,-37\}

Sektor 4 alpha= 12 beta= 4 k= \{6,-40,-37\}

Sektor 5 alpha= 13 beta= 3 k= \{-46,-40,-37\}

Sektor 6 alpha= 9 beta= 5 k= \{-34,-28,-37\}

m= 4

Sektor 1 alpha= 20 beta= 3 k= \{-10,-44,-19\}

Sektor 2 alpha= 26 beta= 9 k= \{-34,-44,-19\}

Sektor 3 alpha= 21 beta= 5 k= \{6,-4,-19\}

Sektor 4 alpha= 14 beta= 4 k= \{6,-28,-19\}

Sektor 5 alpha= 15 beta= 3 k= \{-34,-28,-19\}

Sektor 6 alpha= 11 beta= 6 k= \{-10,-4,-19\}

m= 5

Sektor 1 alpha= 23 beta= 3 k= \{26,16,11\}

Sektor 2 alpha= 31 beta= 11 k= \{-10,16,11\}

Sektor 3 alpha= 23 beta= 5 k= \{6,32,11\}

Sektor 4 alpha= 19 beta= 5 k= \{6,-4,11\}

Sektor 5 alpha= 17 beta= 3 k= \{-10,-4,11\}

Sektor 6 alpha= 13 beta= 7 k= \{26,32,11\}\\

For u={[}m,b,a,b,a,...{]} and m=1...a with the following 

PARAMETERS: a= 3 b= 5 k= 1 (k=1 means smoothness is $C^{1}$)\\

m= 1

Sektor 1 alpha= 50 beta= 3 k= \{-26,-52,-21\}

Sektor 2 alpha= 11 beta= 4 k= \{-26,-52,-21\}

Sektor 3 alpha= 11 beta= 4 k= \{10,-16,-21\}

Sektor 4 alpha= 38 beta= 15 k= \{10,-16,-21\}

Sektor 5 alpha= 37 beta= 3 k= \{-26,-16,-21\}

Sektor 6 alpha= 47 beta= 21 k= \{-26,-16,-21\}

m= 2

Sektor 1 alpha= 16 beta= 3 k= \{-6,-32,-11\}

Sektor 2 alpha= 17 beta= 6 k= \{-26,-32,-11\}

Sektor 3 alpha= 12 beta= 4 k= \{10,4,-11\}

Sektor 4 alpha= 10 beta= 4 k= \{10,-16,-11\}

Sektor 5 alpha= 12 beta= 3 k= \{-26,-16,-11\}

Sektor 6 alpha= 12 beta= 6 k= \{-6,4,-11\}

m= 3

Sektor 1 alpha= 16 beta= 3 k= \{34,28,19\}

Sektor 2 alpha= 20 beta= 7 k= \{-6,28,19\}

Sektor 3 alpha= 14 beta= 4 k= \{10,44,19\}

Sektor 4 alpha= 11 beta= 4 k= \{10,4,19\}

Sektor 5 alpha= 13 beta= 3 k= \{-6,4,19\}

Sektor 6 alpha= 9 beta= 5 k= \{34,44,19\}

\newpage

\subsection{CLLP for the Classic 18-Cycle} \label{CLLP-18cycle}

The classic 18-cycle has the sector sequence  4343-425-43434343-425, i.e. we start in sector 4 with u={[}3,5,3,5,...{]} which is around 3.18819:\\

Sector 4, u= 3.18819

Passage A: GP from Sector 4 to 3. In sector 4, the eigenvalues are:

lc= 1.5418 lt= 1.91557 ln= -1.33278 la= 0.209019

Passage A: GP from Sector 3 to 4. In sector 3, the eigenvalues are:

lc= -2.02211 lt= 3.44689 ln= 2.39822 la= 0.37611

The Eigenvalues are \{1.4570,1.0000\}

The Eigenvectors are \{\{0.167325,1.0000\},\{0,1.0000\}\}\\

Sector 4, u= 2.18819

Passage A: GP from Sector 4 to 3. In sector 4, the eigenvalues are:

lc= 2.02211 lt= 1.42478 ln= -1.646 la= 0.37611

Passage A: GP from Sector 3 to 4. In sector 3, the eigenvalues are:

lc= -2.81366 lt= 3.15683 ln= 3.64699 la= 0.833333

The Eigenvalues are \{1.8416,1.0000\}

The Eigenvectors are \{\{0.62500,1.0000\},\{0,1.0000\}\}\\

Sector 4, u= 1.18819

Passage B2: GP from Sector 4 to 2.In sector 4,the eigenvalues are:

lc= 2.81366 lt= 0.343171 ln= -1.98032 la= 0.833333

Passage B2: GP from Sector 2 to 5.In sector 2,the eigenvalues are:

lc= -3.37457 lt= 1.00965 ln= 5.82633 la= 2.45176

Passage B2: tGP from Sector 5 to 4.In sector 5,the eigenvalues are:

lc= 3.37457 lt= -2.36492 ln= -0.922814 la= 2.45176

The Eigenvalues are \{3.1761,-1.9879\}

The Eigenvectors are \{\{1.01227,1.0000\},\{-0.239459,1.0000\}\}\\

Sector 4, u= 5.31366

Passage A: GP from Sector 4 to 3. In sector 4, the eigenvalues are:

lc= 1.00965 lt= 2.36492 ln= -0.922814 la= 0.0868342

Passage A: GP from Sector 3 to 4. In sector 3, the eigenvalues are:

lc= -1.20737 lt= 3.41557 ln= 1.33278 la= 0.125411

The Eigenvalues are \{1.2318,1.0000\}

The Eigenvectors are \{\{0.045618,1.0000\},\{0,1.0000\}\}\\

Sector 4, u= 4.31366

Passage A: GP from Sector 4 to 3. In sector 4, the eigenvalues are:

lc= 1.20737 lt= 2.2082 ln= -1.08196 la= 0.125411

Passage A: GP from Sector 3 to 4. In sector 3, the eigenvalues are:

lc= -1.49614 lt= 3.45384 ln= 1.6923 la= 0.196156

The Eigenvalues are \{1.3018,1.0000\}

The Eigenvectors are \{\{0.075222,1.0000\},\{0,1.0000\}\}\\

Sector 4, u= 3.31366

Passage A: GP from Sector 4 to 3. In sector 4, the eigenvalues are:

lc= 1.49614 lt= 1.9577 ln= -1.29998 la= 0.196156

Passage A: GP from Sector 3 to 4. In sector 3, the eigenvalues are:

lc= -1.94792 lt= 3.45473 ln= 2.29407 la= 0.346154

The Eigenvalues are \{1.4322,1.0000\}

The Eigenvectors are \{\{0.150000,1.0000\},\{0,1.0000\}\}\\

Sector 4, u= 2.31366

Passage A: GP from Sector 4 to 3. In sector 4, the eigenvalues are:

lc= 1.94792 lt= 1.50681 ln= -1.60176 la= 0.346154

Passage A: GP from Sector 3 to 4. In sector 3, the eigenvalues are:

lc= -2.69398 lt= 3.23295 ln= 3.43668 la= 0.742693

The Eigenvalues are \{1.7612,1.0000\}

The Eigenvectors are \{\{0.49035,1.0000\},\{0,1.0000\}\}\\

Sector 4, u= 1.31366

Passage B2: GP from Sector 4 to 2.In sector 4, the eigenvalues are:

lc= 2.69398 lt= 0.538969 ln= -1.95129 la= 0.742693

Passage B2: GP from Sector 2 to 5.In sector 2, the eigenvalues are:

lc= -3.45737 lt= 1.5418 ln= 5.58196 la= 2.12459

Passage B2: GP from Sector 5 to 4.In sector 5, the eigenvalues are:

lc= 3.45737 lt= -1.91557 ln= -1.33278 la= 2.12459

The Eigenvalues are \{2.8062,-1.4925\}

The Eigenvectors are \{\{2.03045,1.0000\},\{-0.262732,1.0000\}\}\\

Sector 4, u= 3.18819

The Eigenvalues are \{514.49,0.55783\}

The Eigenvectors are \{\{0.448591,1.0000\},\{-0.414926,1.0000\}\}


\newpage

\begin{thebibliography}{99}








\bibitem{ode-reference} H. Amann.
\newblock {\em Ordinary differential equations. }
\newblock Walter de Gruyter, 1990.


\bibitem{beg10} F.~B\'eguin.
\newblock {\em Aperiodic oscillatory asymptotic behavior for some Bianchi
spacetimes. }
\newblock Class.\ Quantum\ Grav.\ 27, 2010.

\bibitem{bkl70} V.A.~Belinski\v{\i}, I.M.~Khalatnikov, and E.M.~Lifshitz.
\newblock {\em Oscillatory approach to a singular point in the relativistic cosmology. }
\newblock Adv.~Phys. 19, 1970.

\bibitem{bkl82} V.A.~Belinski\v{\i}, I.M.~Khalatnikov, and E.M.~Lifshitz.
\newblock {\em A general solution of the Einstein equations with a time singularity. }
\newblock Adv.~Phys.~ 31, 1982.

\bibitem{bianchi-orginal} L. Bianchi.
\newblock {\em Sugli spazii a tre dimensioni che ammettono un gruppo continuo di movimenti. (On the spaces of three dimensions that admit a continuous group of movements). }
\newblock Soc. Ital. Sci. Mem. di Mat. 11, 267, 1898.

\bibitem{bron94}  I. U. Bronstein, A. Y. Kopanskii. 
\newblock {\em Smooth Invariant Manifolds and Normal Forms. }
\newblock World Scientific, 1994.

\bibitem{buchner-evolution} J. Buchner.
\newblock {\em The simplest form of Evolution Equations containing both Gowdy and the exceptional Bianchi cosmological models. }
\newblock http://dynamics.mi.fu-berlin.de/preprints/buchner-g2-equations.pdf, 2010.


\bibitem{einstein15a} A. Einstein. 
\newblock {\em Die Feldgleichungen der Gravitation. }
\newblock Sitzungsberichte der Königlich Preußischen Akademie der Wissenschaften (Berlin), p. 844-847, 1915.

\bibitem{einstein15b} A. Einstein.
\newblock {\em  Zur allgemeinen Relativitätstheorie. }
\newblock Sitzungsberichte der Königlich Preußischen Akademie der Wissenschaften (Berlin), p. 778-786, 1915.

\bibitem{einstein16} A. Einstein.
\newblock {\em Die Grundlage der Allgemeinen Relativitätstheorie. }
\newblock Annalen der Physik, Volume 354, Issue 7, p. 769–822, 1916.

\bibitem{elsetal02} H.~van~Elst, C.~Uggla, and J.~Wainwright.
\newblock {\em Dynamical systems approach to G2 cosmology. }
\newblock Class.\ Quantum Grav. 19, 2002. 



\bibitem{lafo04} S. Gallot, D. Hulin, J. Lafontaine.
\newblock {\em Riemannian Geometry. }
\newblock Springer, 3rd edition, 2004.

\bibitem{grob59} D. M. Grobman.
\newblock {\em Homeomorphism of systems of differential equations. }
\newblock Doklady Akad. Nauk SSSR 128, p. 880–881, 1959.

\bibitem{grob62} D. M. Grobman.
\newblock {\em Topological classification of neighborhoods of a singularity in n-space. }
\newblock Mat. Sb. (N.S.), 56(98):1, p. 77–94, 1962.


\bibitem{hart60b} P. Hartman.
\newblock {\em A lemma in the theory of structural stability of differential equations. }
\newblock Proc. A.M.S. 11 (4): 610–620, 1960.

\bibitem{hartman60a} P. Hartman.
\newblock {\em On local homeomorphisms of Euclidean spaces. }
\newblock Bol. Soc. Mat. Mexicana 5, 220-241, 1960.

\bibitem{heiugg09b} J.~M.~Heinzle and C.~Uggla.
\newblock {\em Mixmaster: Fact and Belief. }
\newblock Class.\ Quantum\ Grav.\ 26, 2009.

\bibitem{heietal09} J.M. Heinzle, C. Uggla, and N. R\"ohr.
\newblock {\em The cosmological billiard attractor. }
\newblock Adv.\ Theor.\ Math.\ Phys. 13, 2009.

\bibitem{heietal12} J.~M.~Heinzle, C.~Uggla, W.~C.~Lim.
\newblock {\em Spike Oscillations. }
\newblock Phys.\ Rev.\ D 86, 2012.

\bibitem{heiugg13} J.~M. Heinzle and C. Uggla.
\newblock {\em Spike statistics. }
\newblock Gen.\ Rel.\ Grav.\ 45, 2013. 

\bibitem{hew03} C. G. Hewitt, J.T. Horwood, J. Wainwright.
\newblock {\em Asymptotic Dynamics of the Exceptional Bianchi Cosmologies. }
\newblock Classical and Quantum Gravity, 20, p. 1743-56, 2003.

\bibitem{hpugh68} Morris W. Hirsch and Charles C. Pugh.
\newblock {\em Stable Manifolds for Hyperbolic Sets. }
\newblock Bull. Amer. Math. Soc. Volume 75, Number 1, p. 149-152, 1969.

\bibitem{hpugh70} Morris W. Hirsch and Charles C. Pugh.
\newblock {\em Stable Manifolds and Hyperbolic Sets. }
\newblock In: Global Analysis, Proceedings of the Symposium, vol. 14, pp. 133-163, AMS, Providende, RI, 1970.

\bibitem{hpugh-shub-book} M. W. Hirsch, C. C. Pugh, and M. Shub.
\newblock {\em Invariant manifolds. }
\newblock Lecture Notes in Mathematics, Vol. 583, Springer-Verlag, Berlin, 1977.


\bibitem{khin49} A. Khintchine.
\newblock {\em Kettenbrüche. }
\newblock B. G. Teubner, Leipzig, 1956 (2nd edition).


\bibitem{liebscher-habil} S. Liebscher.
\newblock {\em Bifurcation without parameters. }
\newblock Habilitationsschrift, Freie Universitiät Berlin, 2012. 

\bibitem{lieetal10} S.~Liebscher, J.~H\"arterich, K.~Webster, and M.~Georgi.
\newblock {\em Ancient Dynamics in Bianchi Models: Approach to Periodic Cycles. }
\newblock Commun.\ Math.\ Phys.\ 305, 2011.

\bibitem{lieetal12} S.~Liebscher, A.~D.~Rendall, and S.~B.~Tchapnda.
\newblock {\em Oscillatory singularities in Bianchi models with magnetic fields. }
\newblock arXiv:1207.2655, 2012.

\bibitem{lk63} E.M.~Lifshitz and I.M.~Khalatnikov.
\newblock {\em Investigations in relativistic cosmology. }
\newblock Adv.~Phys. 12, 1963.

\bibitem{lim04} W.C.~Lim.
\newblock {\em The Dynamics of Inhomogeneous Cosmologies. }
\newblock Ph.~D. thesis, University of Waterloo, 2004; arXiv:gr-qc/0410126.

\bibitem{lim08} W.C.~Lim.
\newblock {\em New explicit spike solution -- non-local component of the generalized Mixmaster  attractor. }
\newblock Class. Quantum Grav. 25, 2008.

\bibitem{limetal09} W.C.~Lim, L.~Andersson, D.~Garfinkle and F.~Pretorius.
\newblock {\em Spikes in the Mixmaster regime of $G_2$ cosmologies. }
\newblock Phys. Rev. D 79, 2009.


\bibitem{matlab-ode} Matlab Documentation Center.\\
\newblock {\em Numerical Integration and Differential Equations. }
\newblock http://www.mathworks.com/help/matlab/ref/ode113.html, 2013.

\bibitem{mis69a} C.~W.~Misner.
\newblock {\em Mixmaster universe. }
\newblock Phys.\ Rev.\ Lett. 22, 1969.

\bibitem{wheeler73} C.W. Misner, K.S. Thorne, and J.A. Wheeler.
\newblock {\em Gravitation}.
\newblock W.H. Freeman and Company, San Francisco, 1973.




\bibitem{takens93} J. Palis, F. Takens. 
\newblock {\em Hyperbolicity and sensitive chaotic dynamics at homoclinic bifurcations. }
\newblock Cambridge University Press, 1993.

\bibitem{per54} O. Perron.
\newblock {\em Die Lehre von den Kettenbrüchen.  }
\newblock B. G. Teubner, Leipzig 1954 (3rd edition).


\bibitem{reitru10} M.~Reiterer and E.~Trubowitz.
\newblock {\em The BKL Conjectures for Spatially Homogeneous
Spacetimes. }
\newblock arXiv:1005.4908v2, 2010.

\bibitem{ren97} A.~D.~Rendall.
\newblock {\em Global dynamics of the mixmaster model. }
\newblock Class. Quantum Grav. 14, 1997. 

\bibitem{ren05} A. D. Rendall.
\newblock {\em The nature of spacetime singularities. }
\newblock In: 100 years of relativity, p. 76 - 92, World Scientific, 2005.

\bibitem{ren08} A.D.~Rendall.
\newblock {\em Partial Differential Equation in General Relativity. }
\newblock Oxford University Press, Oxford, 2008.

\bibitem{rin01} H. Ringstr\"om.
\newblock {\em The Bianchi IX attractor. }
\newblock Annales Henri Poincar\'e 2, 2001.

\bibitem{rin09} H.~Ringstr\"om.
\newblock {\em The Cauchy Problem in General Relativity. }
\newblock ESI Lectures in Mathematics and Physics, 2009.

\bibitem{rob95} C. Robinson. 
\newblock {\em Dynamical Systems: Stability, Symbolic Dynamics, and Chaos. }
\newblock CRC Press, 1995.


\bibitem{sell85} G R Sell. 
\newblock {\em Obstacles to Linearization. }
\newblock Differential Equations, Vol 20, pages 341-345, 1985.

\bibitem{SSTC98} L.~P. Shilnikov, A.~L. Shilnikov, D.~V. Turaev, and L.~O. Chua.
\newblock {\em Methods of Qualitative Theory in Nonlinear Dynamics I.}
\newblock  Volume~4 of Series on Nonlinear Science, Series A, World Scientific, 1998.

\bibitem{shosh72} A. Shoshitaishvili.
\newblock {\em Bifurcations of topological type at singular points of parametrized vector fields. }
\newblock Func anal Appl 6:169-170, 1972.

\bibitem{shosh75} A. Shoshitaishvili.
\newblock {\em Bifurcations of topological type of a vector field near a singular point. }
\newblock Trudy Petrovsky seminar, vol. 1, Moscow University Press, Moscow, pp. 279-309, 1975.

\bibitem{shub-book} M.Shub.
\newblock {\em Global Stability of Dynamical Systems. }
\newblock Springer, 1986.

\bibitem{spivak99} M. Spivak.
\newblock {\em A Comprehensive Introduction to Differential Geometry. }
\newblock Publish or Perish; 3rd edition, 1999.

\bibitem{stern57} S. Sternberg.
\newblock {\em Local Contractions and a Theorem of Poincare. }
\newblock American Journal of Mathematics, Vol. 79, No. 4, pp. 809-824, 1957.

\bibitem{stern58} S. Sternberg.
\newblock {\em On the Structure of Local Homeomorphisms of Euclidean n-Space. }
\newblock American Journal of Mathematics, Vol. 80, No. 3, pp. 623-631, 1958.


\bibitem{takens71} F. Takens. 
\newblock {\em Partially Hyperbolic Fixed Points. }
\newblock Topology Vol.10, 1971. 


\bibitem{ugg13a}  C.~Uggla.
\newblock {\em Spacetime singularities: Recent developments. }
\newblock Int.\ J.\ Mod.\ Phys.\ D\ 22, 2013.

\bibitem{ugg13b}  C.~Uggla.
\newblock {\em Recent developments concerning generic spacelike singularities. }
\newblock Plenary Contribution to ERE2012, http://arxiv.org/abs/1304.6905, 2013.

\bibitem{uggetal03} C.~Uggla,, H.~van Elst, J.~Wainwright and G.F.R.~Ellis.
\newblock {\em The past attractor in inhomogeneous cosmology. }
\newblock Phys.\ Rev.\ D 68, 2003.


\bibitem{vanderbauwhede} Vanderbauwhede, A. 
\newblock {\em Centre Manifolds, Normal Forms and Elementary Bifurcations. }
\newblock Dynamics Reported, 2, 89-169, 1989.


\bibitem{waiell97} J. Wainwright and G.F.R. Ellis.
\newblock {\em Dynamical systems in cosmology. }
\newblock Cambridge University Press, Cambridge, 1997.

\bibitem{waihsu89} J. Wainwright and L. Hsu.
\newblock {\em A dynamical systems approach to Bianchi cosmologies:
orthogonal models of class A. }
\newblock Class.\ Quantum Grav. 6, 1989. 

\bibitem{wald84} R. M. Wald.
\newblock {\em General Relativity. }
\newblock University Of Chicago Press, 1984.

\bibitem{wheeler10} J.A. Wheeler.
\newblock {\em Geons, Black Holes, and Quantum Foam: A Life in Physics. }
\newblock W. W. Norton \& Company, 2010.



\end{thebibliography}
\end{document}